\RequirePackage{ifpdf}
\ifpdf 
\documentclass[pdftex]{sigma}
\else
\documentclass{sigma}
\fi

\DeclareMathOperator{\ad}{ad}

\DeclareMathOperator{\diag}{diag}
\DeclareMathOperator{\Ad}{Ad}
\DeclareMathOperator{\End}{End}
\DeclareMathOperator{\pr}{pr}
\DeclareMathOperator{\spann}{span}

\begin{document}

\allowdisplaybreaks

\renewcommand{\thefootnote}{$\star$}

\renewcommand{\PaperNumber}{093}

\FirstPageHeading

\ShortArticleName{Compact  Riemannian Manifolds  with Homogeneous Geodesics}

\ArticleName{Compact  Riemannian Manifolds\\  with Homogeneous Geodesics\footnote{This paper is a
contribution to the Special Issue ``\'Elie Cartan and Dif\/ferential Geometry''. The
full collection is available at
\href{http://www.emis.de/journals/SIGMA/Cartan.html}{http://www.emis.de/journals/SIGMA/Cartan.html}}}

\Author{Dmitri\u\i\ V. ALEKSEEVSKY~$^\dag$ and Yuri\u\i\ G. NIKONOROV~$^\ddag$}

\AuthorNameForHeading{D.V.~Alekseevsky and Yu.G.~Nikonorov}

\Address{$^\dag$~School of Mathematics and Maxwell Institute for
Mathematical Studies, Edinburgh University,\\
\hphantom{$^\dag$}~Edinburgh EH9 3JZ, United Kingdom}
\EmailD{\href{mailto:D.Aleksee@ed.ac.uk}{D.Aleksee@ed.ac.uk}}

\Address{$^\ddag$~Volgodonsk Institute of Service (branch) of South Russian State
University of Economics\\
\hphantom{$^\ddag$}~and Service,  16 Mira Ave., Volgodonsk, Rostov region, 347386, Russia}
\EmailD{\href{mailto:nikonorov2006@mail.ru}{nikonorov2006@mail.ru}}

\ArticleDates{Received April 22, 2009, in f\/inal form September 20, 2009;  Published online September 30, 2009}

\Abstract{A homogeneous Riemannian space $(M= G/H, g)$ is called a
geodesic orbit  space (shortly, GO-space) if any geodesic is  an
orbit of one-parameter subgroup of the isometry group $G$.  We
study   the structure of  compact GO-spaces  and give  some
suf\/f\/icient conditions for existence and non-existence of an
invariant metric $g$  with  homogeneous geodesics on a~homogeneous
space of a compact Lie group $G$. We give  a classif\/ication of
compact  simply connected GO-spaces $(M = G/H,g)$ of positive
Euler characteristic. If   the group $G$ is simple  and  the
metric $g$ does not come  from  a bi-invariant metric of $G$, then
$M$  is one of the f\/lag  manifolds $M_1=SO(2n+1)/U(n)$ or $M_2=
Sp(n)/U(1)\cdot Sp(n-1)$  and $g$ is any invariant metric on $M$
which depends on two  real parameters. In both cases,  there
exists unique (up to a scaling) symmetric metric $g_0$ such that
$(M,g_0)$ is  the  symmetric space $M = SO(2n+2)/U(n+1)$ or,
respectively, $\mathbb{C}P^{2n-1}$. The manifolds $M_1$, $M_2$ are
weakly symmetric spaces.}


\Classification{53C20; 53C25; 53C35}

\renewcommand{\thefootnote}{\arabic{footnote}}
\setcounter{footnote}{0}

\section{Introduction}\label{section0}

A Riemannian manifold $(M,g)$ is called a  manifold  with homogeneous
geodesics or geodesic orbit manifold (shortly, GO-manifold) if all its geodesic are orbits  of
one-parameter  groups of isometries  of $(M,g)$.  Such manifold is a  homogeneous manifold  and  can  be identif\/ied
  with a~coset space $M = G/H$  of a transitive Lie group $G$ of isometries.
A Riemannian  homogeneous space  $(M = G/H, g^M)$ of a group $G$ is  called a  space with  homogeneous geodesics
(or  geodesic orbit  space, shortly,  GO-space)
if  any geodesic is an orbit
of a one-parameter subgroup of  the group~$G$.
This terminology was
introduced by O.~Kowalski and L.~Vanhecke in~\cite{KV}, who
initiated a~systematic study of such spaces.

Recall  that  homogeneous geodesics correspond to ``relative equilibria''  of  the  geodesic f\/low,
considered as a hamiltonian system on
the  cotangent bundle. Due to this, GO-manifolds can be  characterized  as   Riemannian manifolds such that all
integral  curves of the geodesic f\/low  are  relative equilibria.

GO-spaces  may be considered as  a natural generalization of   symmetric  spaces, classif\/ied by \'E.~Cartan~\cite{Ca}.
Indeed, a simply connected  symmetric space can be def\/ined as a Riemannian manifold $(M,g)$ such  that any geodesic
$\gamma \subset M$ is  an orbit of  one-parameter  group $g_t$ of  transvections, that is one-parameter  group of
isometries which preserves $\gamma$ and induces the  pa\-ral\-lel  transport along $\gamma$.
If we remove the assumption that $g_t$ induces the pa\-ral\-lel transport, we  get the notion of a GO-space.

The class of GO-spaces is much larger  then the  class of symmetric spaces. Any homogeneous space
$M = G/H$ of a compact Lie group $G$  admits a metric $g^M$ such that $(M,g^M)$ is a GO-space.
It is suf\/f\/icient to take  the metric~$g^M$ which is induced with a bi-invariant Riemannian metric~$g$ on the Lie group~$G$ such that
$ (G,g) \to (M=G/H, g^M)$ is a~Riemannian submersion
with totally geodesic f\/ibres. Such  GO-space $(M = G/H, g^M)$ is called  a {\bf normal homogeneous space}.

More generally,   any  naturally  reductive  manifold  is a  geodesic orbit manifold.
Recall that a~Riemannian manifold $(M,g^M)$ is called
{\bf naturally reductive}  if it admits a transitive Lie group~$G$ of isometries with a bi-invariant
pseudo-Riemannian metric $g$, which induces the metric~$g^M$ on $M = G/H$, see \cite{KN,Bes}.
The f\/irst example of non naturally reductive GO-manifold had been constructed by  A.~Kaplan~\cite{Kap}.
An important class of GO-spaces   consists  of weakly symmetric spaces,  introduced by A.~Selberg~\cite{S}.
A homogeneous Riemannian  space  $(M = G/H, g^M)$ is a {\bf  weakly symmetric space}
if any two points $p,q \in M$ can be interchanged by
an isometry $ a \in G$. This  property does not depend on the particular invariant metric $g^M$.
Weakly symmetric spaces $M= G/H$
have many interesting properties
(for example, the  algebra  of $G$-invariant dif\/ferential operators on $M$ is commutative,  the  representation of
$G$ in the space $L^2(M) $ of function is multiplicity free, the algebra  of $G$-invariant Hamiltonians on $T^*M$
with respect to Poisson bracket is commutative)  and  are closely related  with spherical spaces, commutative
spaces and Gelfand pairs etc., see  the book by J.A.~Wolf~\cite{W1}. The classif\/ication of weakly symmetric
reductive homogeneous  spaces was  given by O.S.~Yakimova~\cite{Yak},
see  also~\cite{W1}.

In \cite{KV}, O.~Kowalski and L.~Vanhecke classif\/ied all GO-spaces
of dimension $\leq 6$.   C.~Gordon~\cite{Gor96} reduced  the
classif\/ication of GO-spaces to the
classif\/ication of GO-metrics on nil\-ma\-ni\-folds, compact GO-spaces  and  non-compact
GO-spaces  of  non-compact semisimple Lie group.
She  described  GO-metrics on nilmanifolds. They exist only  on two-step nilponent
nilmanifolds. She also presented some constructions of GO-metrics on homogeneous compact
manifolds  and non compact manifolds of a semisimple  group.

Many  interesting results  about  GO-spaces one can f\/ind  in \cite{BKV,DuKoNi,Zi96,Ta,Tam}, where there are  also extensive references.

Natural generalizations of normal homogeneous Riemannian manifolds are $\delta$-homoge\-neous Riemannian manifolds,
studied in \cite{BerNik, BerNik3, BerNikNik}.
Note that the class of $\delta$-homoge\-neous Riemannian manifolds is a proper subclass
of the class of geodesic orbit spaces with non-negative sectional curvature (see the quoted papers for
further properties of $\delta$-homogeneous Riemannian manifolds).

In \cite{AA}, a   classif\/ication  of  non-normal invariant GO-metrics on  f\/lag manifolds $M =G/H$ was  given.
The  problem reduces to the case when the  (compact)  group $G$ is simple. There  exist only two
series of  f\/lag manifolds of a simple group which admit  such metric, namely  weakly  symmetric  spaces
$M_1= SO(2n+1)/U(n)$  and $M_2= Sp(n)/U(1) \cdot Sp(n-1)$, equipped  with  any (non-normal)
invariant metric (which depends  on two real parameters).
Moreover,  there exists  unique  (up to a scaling) invariant metric $g_0$,   such that the Riemannian
manifolds $(M_i, g_0)$  are isometric to the  symmetric spaces $SO(2n+2)/U(n+1)$ and
$\mathbf{C}P^{2n-1}  = SU(2n)/U(2n-1)$, respectively.

The main  goal of  this paper is a generalization of this result to the  case of compact homogeneous manifolds
of positive Euler  characteristic. We prove that the  weakly symmetric  manifolds  $M_1$, $M_2$   exhaust
all simply connected compact irreducible Riemannian  non-normal GO-manifolds of positive Euler
characteristic.

We  indicate now the idea  of the proof.   Let $(M = G/H, g^M)$ be a  compact  irreducible non-normal
GO-space of positive Euler characteristic.
Then  the stability subgroup~$H$ has maximal rank, which implies  that $G$ is  simple.
We prove that there is rank~2 regular simple subgroup
$G'$ of $G$  (associated with a rank~2  subsystem $R'$  of the root system~$R$  of
the Lie algebra  $\mathfrak{g} = {\rm Lie}(G)$)  such that
the  orbit $M' = G'o = G'/H'$ of the point $o = eH \in M$  (with the induced metric) is  a~non-normal GO-manifold.
Using \cite{AA,BerNik}, we prove that
the only such manifold $M'$  is $SU(5)/U(2)$. This  implies that the  root system $R$ is not
simply-laced  and  admits a  ``special'' decomposition
$R = R_0 \cup R_1 \cup R_2$  into a  disjoint union  of three  subsets, which satisf\/ies  some properties.
We determine  all such special decompositions
of  irreducible   root systems and  show that  only  root systems of type  $B_n$ and $C_n$  admit
special decomposition   and associated  homogeneous
manifolds  are  $M_1$  and $M_2$.

The structure of the paper is the following.  We f\/ix notations and recall  basic def\/initions in Section~\ref{section1}.
Some  standard facts about totally geodesic submanifolds of a homogeneous
Riemannian spaces  are collected in Section~\ref{section2}.
We  discuss  some properties of compact  GO-spaces  in Section~\ref{section3}.
These results  are used in Section~\ref{section4} to derive suf\/f\/icient  conditions  for  existence and non-existence
of  a non-normal GO-metric on
a homogeneous manifold of a compact group.
Section~\ref{section5} is devoted to classif\/ication of compact GO-spaces with positive Euler characteristic.

\section{Preliminaries and notations}\label{section1}

Let $M= G/H$ be a homogeneous space of a compact connected Lie group $G$.
We will denote by $b= \langle \cdot ,\cdot \rangle$ a f\/ixed  $\Ad_G$-invariant Euclidean metric
 on the Lie algebra
$\mathfrak{g}$ of $G$ (for example, the minus Killing form if $G$ is semisimple) and by
\begin{gather}\label{reductivedecomposition}
\mathfrak{g}=\mathfrak{h} + \mathfrak{m}
\end{gather}
the  associated $b$-orthogonal  reductive decomposition, where
$\mathfrak{h} = {\rm Lie}(H) $. An invariant Riemannian metric $g^M$
on $M$  is determined by an $\Ad_H$-invariant Euclidean metric $g =
(\cdot,\cdot)$ on  the space $\mathfrak{m}$ which is identif\/ied with
the tangent space $T_oM$ at the initial point $o = eH$.

If $\mathfrak{p}$ is a subspace of $\mathfrak{m}$, we will denote by
$X_{\mathfrak{p}}$ the $b$-orthogonal projection of a vector $X \in
\mathfrak{g}$ onto $\mathfrak{p}$, by $b_{\mathfrak{p}}$ the
restriction of the symmetric bilinear form to $\mathfrak{p}$ and by
$A^{\mathfrak{p}} =\pr_{\mathfrak{p}} \circ A \circ
\pr_{\mathfrak{p}}$ the projection of  an endomorphism $A$ to
$\mathfrak{p}$. If $g$ is a $\Ad_H$-invariant metric, the quotient
\[
A ={ b_{\mathfrak{m}}^{-1}} \circ g
\]
 is  an
$\Ad_H$-equivariant  symmetric positively def\/ined  endomorphism  on
$\mathfrak{m}$, which we call the { \bf metric endomorphism}.
Conversely, any  such equivariant positively def\/ined  endomorphism~$A$
of $\mathfrak{m}$ def\/ines an invariant metric $g = b \circ A = b (A
\cdot, \cdot)$ on $\mathfrak{m}$, hence an invariant Riemannian
metric~$g^M$ on~$M$.

\begin{lemma}\label{egen}
Let $(M=G/H, g^M)$ be a compact homogeneous Riemannian space with  metric endomorphism $A$  and
\begin{equation}\label{eigensum}
\mathfrak{m}=\mathfrak{m}_1 \oplus \mathfrak{m}_2 \oplus \cdots \oplus \mathfrak{m}_k,
\end{equation}
 the $A$-eigenspace   decomposition  such that  $A|_{\mathfrak{m}_i} = \lambda_i \cdot{\bf 1}_{\mathfrak{m}_i}$.
Then
\begin{equation}\label{orthogonality}
(\mathfrak{m}_i,\mathfrak{m}_j)=\langle \mathfrak{m}_i,\mathfrak{m}_j \rangle=0
\end{equation}
and $\Ad_H$-modules $\mathfrak{m}_i$ satisfy
$[\mathfrak{m}_i,\mathfrak{m}_j]\subset \mathfrak{m}$ for $i\neq j$.
\end{lemma}

\begin{proof}
 Since $A$ commute with $\Ad_H$,   eigenspaces $\mathfrak{m}_i$ are $\Ad_H$-invariants  and  for
 $X\in \mathfrak{m}_i$, $Y\in \mathfrak{m}_j$, $i\neq j$, we get
\[
\lambda_i \langle X,Y\rangle =\langle AX,Y\rangle=(X,Y)=\langle X,AY\rangle=\lambda_j \langle X,Y\rangle.
\]
 This  implies  (\ref{orthogonality}).
The inclusion $[\mathfrak{m}_i,\mathfrak{m}_j]\subset \mathfrak{m}$ follows from the fact that $\mathfrak{m}_j$
is $\Ad_H$-invariant and
$
\langle [\mathfrak{m}_i,\mathfrak{m}_j], \mathfrak{h}\rangle=
\langle \mathfrak{m}_i,[\mathfrak{m}_j, \mathfrak{h}]\rangle=0
$.
\end{proof}

For any subspace $\mathfrak{p} \subset \mathfrak{m}$ we will denote by
$\mathfrak{p}^\perp$ its orthogonal complement with respect to the metric $g$ and by
$\bold{1}_\mathfrak{p}$ the identity operator on $\mathfrak{p}$.

Recall  that  $\Ad_H$-submodules $\mathfrak{p}$, $\mathfrak{q}$ are called { \bf disjoint}
 if they have no non-zero equivalent  submodules.
If $\Ad_H$-module $\mathfrak{m}$ is decomposed into a direct sum
\[
\mathfrak{m}= \mathfrak{m_1}+ \cdots + \mathfrak{m}_k
\]
of  disjoint submodules, then any $\Ad_H$-invariant
metric $g$ and associated metric endomorphism~$A$ have the form
\[
g = g_{\mathfrak{m}_1} \oplus \cdots \oplus  g_{\mathfrak{m}_k} ,
\qquad A = A^{\mathfrak{m}_1} \oplus \cdots \oplus
A^{\mathfrak{m}_k}.
\]

Let $(M= G/H,g^M)$ be a compact homogeneous Riemannian space
with the reductive decomposition (\ref{reductivedecomposition})
and  metric endomorphism $A \in \End({\mathfrak m})$.

We identify  elements $X,Y \in  \mathfrak{g}$ with Killing vector f\/ields on $M$.
Then the covariant derivative~$\nabla_XY$ at the point $o = eH$ is
given by
\begin{equation}\label{connec}
\nabla_XY(o)=-\tfrac{1}{2}[X,Y]_{\mathfrak{m}}+U(X_{\mathfrak{m}},Y_{\mathfrak{m}}),
\end{equation}
where the bilinear symmetric map $U:\mathfrak{m}\times \mathfrak{m}
\rightarrow \mathfrak{m}$ is given by
\begin{equation}\label{connec1}
2(U(X,Y),Z)=( \ad_Z^{\mathfrak{m}}X,Y)+(X,\ad_Z^{\mathfrak{m}}Y)
\end{equation}
for any $X,Y, Z\in \mathfrak{m}$ and $X_{\mathfrak{m}}$ is the
$\mathfrak{m}$-part of a vector $X \in \mathfrak{g} $ \cite{Bes}.

\begin{definition}\label{GOsp}
A homogeneous Riemannian space  $(M=G/H,g^M)$ is called a  space
with homogeneous geodesics  { shortly,  \bf GO-space} if any
geodesic $\gamma $ of $M$ is an orbit of  1-parameter subgroup of
$G$. The invariant  metric $g^M$ is called { \bf GO-metric}.
\end{definition}
   If $G$ is the  full isometry group,   then GO-space is
  called    a {\bf manifold with  homogeneous geodesics} or {\bf GO-manifold}.

\begin{definition}\label{stGOsp}
 A GO-space $(M=G/H,g^M)$  of a simple  compact Lie group $G$ is called a {\bf
proper GO-space} if the metric $g^M$ is not $G$-normal, i.e.\ the
metric endomorphism  $A$ is not a scalar operator.
\end{definition}

\begin{lemma}[\cite{AA}]\label{GO-criterion}
A  compact homogeneous Riemannian space   $(M=G/H,g^M)$ with the
reductive decomposition  \eqref{reductivedecomposition} and  metric
endomorphism  $A$ is GO-space if and
only if  for any $X \in \mathfrak{m}$ there is
$H_X \in \mathfrak{h}$ such
that  one of the following  equivalent conditions holds:
\begin{enumerate}\itemsep=0pt
\item[$i)$] $[H_X +X, A(X)]\in \mathfrak{h}$;

\item[$ii)$]  $([H_X +X,Y]_{\mathfrak{m}},X) =0$ for all $Y\in \mathfrak{m}$.
\end{enumerate}
\end{lemma}

This lemma  shows that  the property to be GO-space  depends only on the reductive
decomposition (\ref{reductivedecomposition})  and the Euclidean metric $g$ on
$\mathfrak{m}$. In other words, if $(M = G/H, g^M)$ is a GO-space, then any locally
isomorphic homogeneous Riemannian space  $(M'= G'/H', g^{M'})$ is a GO-space.  Also  a
direct product of Riemannian manifolds is a manifold with homogeneous geodesics if and
only if each factor is a manifold with homogeneous geodesics.

\section{Totally geodesic orbits in a homogeneous Riemannian space}\label{section2}

In this section we deal with totally geodesic submanifolds of
compact homogeneous Riemannian spaces. This is  a useful tool for
study of  GO-spaces due to the following

\begin{proposition}[\protect{\cite[Theorem 11]{BerNik}}]\label{tot}
Every closed totally geodesic submanifold of a  Riemannian manifold
with  homogeneous geodesics is a manifold with  homogeneous
geodesics.
\end{proposition}

Let $(M = G/H,g^M)$ be a compact Riemannian homogeneous space
with  the reductive decomposition (\ref{reductivedecomposition}).

\begin{definition} A subspace $\mathfrak{p} \subset \mathfrak{m}$
is called { \bf totally geodesic} if it is the tangent space  at $o$
of a~totally geodesic orbit $Ko \subset G/H =M$ of a subgroup $K
\subset G$.
\end{definition}

\begin{proposition}\label{totallygeodesicProp}
A subspace $\mathfrak{p} \subset \mathfrak{m}$
is totally geodesic if  and only if the following two conditions hold:
\begin{enumerate}\itemsep=0pt
\item[$a)$] $\mathfrak{p}$ generates a  subalgebra  of the form $\mathfrak{k} =
\mathfrak{h}' + \mathfrak{p}$, where $\mathfrak{h}'$ is a
subalgebra of $\mathfrak{h}$;

\item[$b)$] the endomorphism $\ad_Z^{\mathfrak{p}} \in \End(\mathfrak{p})$
for  $Z \in \mathfrak{p}^{\perp}$ is $g$-skew-symmetric or,
equivalently,
\[
U(\mathfrak{p}, \mathfrak{p}) \subset \mathfrak{p}.
\]
\end{enumerate}
\end{proposition}

\begin{proof}
If $\mathfrak{p}$ is the  tangent space of the orbit  $Ko = K/H'$, then
 ${\rm Lie}(K)=\mathfrak{k}=\mathfrak{h}'+\mathfrak{p}$, where
$\mathfrak{h}' = {\rm Lie}( H')$ is a subalgebra of $\mathfrak{h}$.
Moreover,  the formulas (\ref{connec}) and (\ref{connec1}) imply
$U(\mathfrak{p}, \mathfrak{p}) \subset \mathfrak{p}$. Conversely,
the conditions $a)$ and $b)$ imply that $\mathfrak{p}$ is the tangent
space  of  the totally geodesic  orbit $Ko$  of the subgroup $K$
generated by the subalgebra~$\mathfrak{k}$.
\end{proof}

\begin{corollary}\qquad\null
\begin{enumerate}
\itemsep=0pt
\item[$i)$] A  subspace  $\mathfrak{p} \subset \mathfrak{m}$ is totally geodesic if a) holds  and $A\mathfrak{p} =
\mathfrak{p}$.

\item[$ii)$] If a totally geodesic  subspace $\mathfrak{p}$ is
$\ad_{\mathfrak{h}}$-invariant  and $A$-invariant, then
\[
 [\mathfrak{h}+ \mathfrak{p}, \mathfrak{p}^{\perp}] \subset \mathfrak{p}^{\perp}.
 \]
\end{enumerate}
\end{corollary}

\begin{proof} $i)$ Assume that  $A\mathfrak{p} = \mathfrak{p}$. Then
$A\mathfrak{p}^{\perp} = \mathfrak{p}^{\perp}$ and $\langle
\mathfrak{p}, \mathfrak{p}^{\perp}\rangle=0$. From $i)$ and
$A\mathfrak{p} = \mathfrak{p}$ we get $\langle Z, [X,AX]\rangle =0$
for any $X\in \mathfrak{p}$ and $Z \in \mathfrak{p}^{\perp}$.
This implies
\begin{gather*}
0=\langle [Z,X], AX\rangle = \langle
[Z,X]_\mathfrak{m},AX\rangle=([Z,X]_\mathfrak{m},X)\\
\phantom{0=\langle [Z,X], AX\rangle }{} = ([Z,X]_\mathfrak{p},X) = (U(X,X),Z).
\end{gather*}

$ii)$ follows from the fact that  the  endomorphisms
$\ad_{\mathfrak{h} + \mathfrak{p}}$   are $b$-skew-symmetric  and
preserves  the  subspace $\mathfrak{p}$. Hence, they preserve its
$b$-orthogonal complement $\mathfrak{p}^{\perp}$.
\end{proof}

\begin{corollary}\label{togc1}
Let $(M = G/H,g)$ be a
 compact Riemannian homogeneous space and
$K$  a connected subgroup of $G$. The orbit $P = Ko = K/H'$ is a
totally geodesic submanifold if and only if the Lie algebra
$\mathfrak{k}$  is consistent with the reductive decomposition
\eqref{reductivedecomposition} $($that is $\mathfrak{k} =
\mathfrak{k}\cap \mathfrak{h}+ \mathfrak{k}\cap \mathfrak{m} =
\mathfrak{h}' + \mathfrak{p})$  and
\[
U(\mathfrak{p}, \mathfrak{p}) \subset \mathfrak{p}
\]
or, equivalently,  the endomorphisms $\ad_Z^{\mathfrak{p}}
\in \End(\mathfrak{p})$, $Z \in \mathfrak{p}^\perp$ are  $g$-skew-symmetric.
\end{corollary}

\section{Properties of  GO-spaces}\label{section3}

\begin{lemma} \label{skew-symmetryLemma}
Let $(M =G/H,g^M)$ be a GO-space with the reductive decomposition
\eqref{reductivedecomposition} and $\mathfrak{m} = \mathfrak{p}+ \mathfrak{q}$
a $g$-orthogonal $\Ad_H$-invariant   decomposition. Then
\[
U(\mathfrak{p},\mathfrak{p} ) \subset \mathfrak{p},\qquad
U(\mathfrak{q}, \mathfrak{q}) \subset \mathfrak{q}
\]
and the endomorphisms
$\ad_{\mathfrak{p}}^{\mathfrak{q}}$, $\ad_{\mathfrak{q}}^{\mathfrak{p}}
$ are skew-symmetric.
\end{lemma}

\begin{proof}  For $X \in \mathfrak{p}$, $Y \in \mathfrak{q} $ we
have
\[
0= \left([Y + H_{Y}, X]_{\mathfrak{m}},Y \right) = - \left(\ad_X Y,Y
\right) =- \left( U(Y,Y), X \right) ,
\]
where $H_Y$ is as in Lemma \ref{GO-criterion}. This shows that
$\ad_{X}^{\mathfrak{q}}$ is skew-symmetric  and
$U(\mathfrak{q},\mathfrak{q} ) \subset \mathfrak{q}$.
\end{proof}

Lemma \ref{skew-symmetryLemma} together with  Proposition \ref{totallygeodesicProp} implies

\begin{proposition}\label{submersion}
Let $(M=G/H, g^M)$ be a GO-space with the reductive decomposition
\eqref{reductivedecomposition}. Then any  connected subgroup $K
\subset G$ which contains $H$ has the totally geodesic orbit $P = Ko
= K/H$ which is GO-space $($with respect to the induced metric$)$.
Moreover, if the space $\mathfrak{p}:= \mathfrak{k} \cap
\mathfrak{m}$ is $A$-invariant, then
\[
 [\mathfrak{k}, \mathfrak{m^{\perp}}] \subset \mathfrak{m}^{\perp}
 \]
and the metric $\bar g : = g|_{\mathfrak{p}^\perp}$ is
$\Ad_K$-invariant and defines an invariant GO-metric $g^N$ on the
homogeneous manifolds $N = G/K$. The projection $\pi : G/H
\rightarrow G/K$ is a Riemannian submersion with totally geodesic
fibers such that the  fibers and the base are GO-spaces.
\end{proposition}

\begin{proof} The f\/irst claim follows from  Lemma~\ref{skew-symmetryLemma}, Lemma
\ref{GO-criterion} and Proposition \ref{totallygeodesicProp}.
 If $A \mathfrak{p} = \mathfrak{p}$, then
$\mathfrak{m }= \mathfrak{p}+ \mathfrak{p}^\perp $ is  a $b$-orthogonal decomposition  and
 since the metric $b$ is $\Ad_G$-invariant, \mbox{$ \Ad_{K}\mathfrak{p}^{\perp} = \mathfrak{p}^\perp$}.
Then Lemma~\ref{skew-symmetryLemma} shows that the metric~$g|_{\mathfrak{p}^\perp}$ is $\Ad_K$-invariant and def\/ines
an invariant metric~$g^N$ on $N = G/K$  such that $N$ becomes GO-space.
\end{proof}

 Note that a subgroup $K \supset H$ is compatible with any invariant metric on $G/H$ if
 $\Ad_H$-modules~$\mathfrak{p}$ and $\mathfrak{m}/ \mathfrak{p}$ are strictly disjoint.
This remark implies

\begin{proposition}  Let $(M=G/H, g)$ be a compact homogeneous Riemannian space.
Then the connected normalizer $N_0(Z)$ of a central subgroup $Z$
of $H$ and the connected normalizer $N_0(H)$ are subgroups
consistent with any invariant metric on~$M$.
\end{proposition}

\begin{proposition} Let $(M=G/H, g^M)$ be a compact GO-space with  metric endomorphism~$A$.
\begin{enumerate}\itemsep=0pt
\item[$i)$] Let $X,Y \in \mathfrak{m}$ be eigenvectors of the metric
endomorphism $A$ with different eigenvalues~$\lambda$,~$\mu$. Then
\[
[X,Y] = \frac{\lambda}{\lambda - \mu}[H,X]+ \frac{\mu}{\lambda -
\mu}[H,Y]
\]
for some $H \in \mathfrak{h}$.

\item[$ii)$] Assume that  the vectors $X$, $Y$ belong to the
$\lambda$-eigenspace $ \mathfrak{m}_{\lambda} $  of $A$  and  $X$ is
 $g$-orthogonal to the subspace $[\mathfrak{h},Y]$. Then
\[
[X,Y] \in
 \mathfrak{h}+ \mathfrak{m}_{\lambda}.
 \]
 \end{enumerate}
\end{proposition}

\begin{proof}
$i)$ Let $X,Y \in  \mathfrak{m} $ be  eigenvectors of $A$ with
dif\/ferent eigenvalues $\lambda$, $\mu$ and $H = H_{X +Y} \in
\mathfrak{h}$ the element def\/ined in Lemma~\ref{GO-criterion}. Then
\begin{gather*}
[H + X + Y, A(X+Y)]
={[H + X + Y, \lambda X + \mu Y]} \\
\phantom{[H + X + Y, A(X+Y)]}{} = \lambda [H,X]+ \mu[H,Y]+ (\mu -
\lambda)[X,Y] \in \mathfrak{h}.
\end{gather*}
By Lemma~\ref{egen}, $[H,X],[H,Y],[X,Y] \in  \mathfrak{m}$  and the right hand side  is zero.

$ii)$ Assume now that  $X,Y \in  \mathfrak{m}_{\lambda}$ satisfy conditions $ii)$
and $Z$ is  an eigenvector of $A$ with an eigenvalue $\mu \neq \lambda$.
Then we have
\begin{gather*}
([X,Y]_{\mathfrak{m}},Z) = \mu \langle[X,Y],Z\rangle =  \mu \langle X,[Y,Z] \rangle =
\frac{\mu}{\lambda}(X,[Y,Z]_{\mathfrak{m}}) \\
\phantom{([X,Y]_{\mathfrak{m}},Z)}{} = \frac{\mu}{\lambda}\left(X, \frac{\lambda}{\lambda - \mu}[H,Y] + \frac{\mu}{\lambda - \mu}[H,Z]\right) = 0.
\end{gather*}
This shows that  $[X,Y] \in \mathfrak{h}+ \mathfrak{m}_{\lambda}$.
\end{proof}

\begin{corollary}\label{CorollaryBrackets} Let $(M=G/H,g^M)$ be a compact GO-space  with  the
reductive decomposition \eqref{reductivedecomposition} and metric endomorphism $A$ and
\begin{equation} \label{A-decomposition}
 \mathfrak{m} = \mathfrak{m}_1 + \cdots  + \mathfrak{m}_k
\end{equation}
 the $A$-eigenspace decomposition such that $A \vert_{\mathfrak{m}_i}
= \lambda_i 1_{\mathfrak{m}_i} $.
Then for  any $\Ad_H$-submodules $\mathfrak{p}_i \subset
\mathfrak{m}_i$, $\mathfrak{p}_j \subset \mathfrak{m}_j$, $i
\neq j$, we have
\[
[\mathfrak{p}_i, \mathfrak{p}_j] \subset  \mathfrak{p}_i + \mathfrak{p}_j.
\]
Moreover, if $\mathfrak{p}$, $\mathfrak{p'}$ are $g$-orthogonal
$\Ad_H$-submodules of $\mathfrak{m}_i$ then
\[
[\mathfrak{p}, \mathfrak{p'}] \subset \mathfrak{h} +  \mathfrak{m}_{i}.
\]
\end{corollary}

\section{Some applications}\label{section4}

\subsection{A suf\/f\/icient condition for non-existence of GO-metric}\label{section4.1}

Here we consider some applications of results of the previous section.
\begin{definition} Let $(M=G/H, g^M)$ be a compact homogeneous Riemannian space.
 A connected closed  Lie subgroup $K \subset G$  which  contains  $H$ is called
 { \bf compatible with the metric $g^M$}
  if the subspace $\mathfrak{p} = \mathfrak{k} \cap \mathfrak{m}$  of $\mathfrak{m}$ is
  invariant under the metric endomorphism $A$.
\end{definition}

 Let $K$, $K'$ be two  subgroups of $G$ which are  compatible  with
 the metric  of a homogeneous Riemannian space $(M = G/H,
 g^M)$. Then we can decompose the space $\mathfrak{m}$ into
a $g$-orthogonal sum of  $A$-invariant $\Ad_H$-modules
 \begin{equation} \label{4termsdecompositron}
 \mathfrak{m} = \mathfrak{q}+ \mathfrak{p}_1 + \mathfrak{p}_2 + \mathfrak{n}
 \end{equation}
where
\begin{gather*}
\mathfrak{q} = \mathfrak{p} \cap  \mathfrak{p}',\qquad
 \mathfrak{p} = \mathfrak{k} \cap \mathfrak{m} = \mathfrak{q}+
\mathfrak{p}_1,\qquad
\mathfrak{p'} = \mathfrak{k'} \cap \mathfrak{m} = \mathfrak{q}+
\mathfrak{p}_2
\end{gather*}
 and $\mathfrak{n}$ is the orthogonal complement to
\[
\mathfrak{p} + \mathfrak{p}' = \mathfrak{q} + \mathfrak{p}_1 +
\mathfrak{p}_2
\]
in $\mathfrak{m}$.

\begin{proposition}\label{2subgroupPropf}
Let $(M=G/H, g^M)$ be a homogeneous Riemannian space,  $K$, $K'$   two
subgroups of $G$  which are compatible with  $g^M$ and \eqref{4termsdecompositron}
 the associated decomposition as above. Then~$\mathfrak{p}_1$,~$\mathfrak{p}_2$,~$\mathfrak{n}$
 are $\Ad_{\widetilde H}$-modules, where $\widetilde H = K \cap K'$ is the
 Lie group with the Lie algebra
 $\widetilde{ \mathfrak{h}} = \mathfrak{h} +
 \mathfrak{q}$,  and
\[
 [\mathfrak{p}_1, \mathfrak{p}_2] \subset \mathfrak{n}.
 \]
Moreover, if $(M=G/H, g^M)$ is a GO-space, then
the  restriction
  $A^{\widetilde{\mathfrak{m}}}$ of the metric endomorphism  to
  $ \widetilde{ \mathfrak{m}} = \mathfrak{p}_1 + \mathfrak{p}_2 +
  \mathfrak{n}$
  commutes with $\Ad_{\widetilde H}|_{\widetilde{\mathfrak{m}}}$  and  for any $\widetilde H$-irreducible
  submodules $\mathfrak{p}'_1 \subset \mathfrak{p}_1$, and
$\mathfrak{p}'_2 \subset \mathfrak{p}_2$ such that
\[
[\mathfrak{p}'_1, \mathfrak{p}'_2]  \neq 0,
\]
the metric endomorphism $A$ is a scalar on the space
\[
\mathfrak{p}'_1 +\mathfrak{p}'_2 + [\mathfrak{p}'_1,
\mathfrak{p}'_2].
\]
\end{proposition}

\begin{proof} Since the decomposition  (\ref{4termsdecompositron}) is $b$-orthogonal,
  we conclude that it is $\Ad_{\widetilde H}$-invariant  and
\begin{gather*}
[\mathfrak{p}, \mathfrak{p}^{\perp}] =
[\mathfrak{p}, \mathfrak{p}_2 + \mathfrak{n}] \subset \mathfrak{p}_2
+ \mathfrak{n},\\
[\mathfrak{p'}, (\mathfrak{p}')^{\perp}] =
[\mathfrak{p}', \mathfrak{p}_1 + \mathfrak{n}] \subset
\mathfrak{p}_1 + \mathfrak{n},
\end{gather*}
by  Proposition
\ref{totallygeodesicProp}. This implies
\[
[\mathfrak{p}_1, \mathfrak{p}_2] \subset \mathfrak{n}.
\]
 If $(M,g^M)$ is a GO-space, then by Proposition \ref{submersion} the metric endomorphism
 $A^{\widetilde {\mathfrak{m}}}$  is $\widetilde H$-invariant. If modules $\mathfrak{p}'_1$, $\mathfrak{p}'_2$ belong to $A$-eigenspaces with dif\/ferent
 eigenvalues, then  by Corollary~\ref{CorollaryBrackets},
\[
[\mathfrak{p}'_1, \mathfrak{p}'_2] \subset  \mathfrak{p}_1 + \mathfrak{p}_2.
\]
 Together with the previous inclusion, it implies  $[\mathfrak{p}'_1, \mathfrak{p}'_2] =0$.
 If  these modules belong to the same eigenspace
 $\mathfrak{m}_{\lambda}$, then by Corollary \ref{CorollaryBrackets},
\[
[\mathfrak{p}'_1, \mathfrak{p}'_2] \subset  \mathfrak{m}_{\lambda}.\tag*{\qed}
\]\renewcommand{\qed}{}
\end{proof}

 As a corollary, we get the following suf\/f\/icient condition that a
 homogeneous manifold $M = G/H$ does not admit a proper
 GO-metric.

 \begin{proposition} \label{2subgroupProp}
 Let $M = G/H$ be a homogeneous space of a compact group $G$ with
 the
 reductive decomposition $\mathfrak{g}  = \mathfrak{h} +
 \mathfrak{m}$. Assume that the Lie algebra $\mathfrak{g}$ has  two subalgebras
$\mathfrak{k} = \mathfrak{h} + \mathfrak{p}$, $\mathfrak{k}'= \mathfrak{h}+ \mathfrak{p'}$ which contain
$\mathfrak{h}$  and  generate $\mathfrak{g}$.
 Let
 \begin{equation}\label{qppn-decomposition}
 \mathfrak{m} =
\mathfrak{q}+ \mathfrak{p}_1 + \mathfrak{p}_2 + \mathfrak{n}, \qquad  \mathfrak{q} = \mathfrak{p}\cap \mathfrak{p'}
\end{equation}
be  the associated $b$-orthogonal decomposition. Assume that there is no commuting
 $\ad_{\mathfrak{h}+ \mathfrak{q}}$ submodules of $\mathfrak{p}_1$  and $\mathfrak{p}_2$.
 Then  for any  GO-metric, defined by an
 operator $A$   which preserves this decomposition,
 $A$ is a scalar operator on $\mathfrak{p}_1 + \mathfrak{p}_2+\mathfrak{n}$.
 In particular, if  $\mathfrak{q}$ is trivial and $\Ad_H$-modules $\mathfrak{p}_1$, $\mathfrak{p}_2$, $\mathfrak{n}$
are strictly non-equivalent, then the only  GO-metric on $M$ is the
normal metric.
 \end{proposition}

\begin{proof}  Let $A$ be  an operator on $\mathfrak{m}$  which
preserves the decomposition (\ref{qppn-decomposition}) and def\/ines a
GO-metric. Then by Proposition \ref{2subgroupPropf},
\[
A|_{\mathfrak{p}_1 + \mathfrak{p}_2 + [\mathfrak{p}_1,
\mathfrak{p}_2]} = \lambda \cdot { \bf 1}
\]
for some $\lambda$. Now
$\mathfrak{p}_1$ and $[\mathfrak{p}_1, \mathfrak{p}_2] \subset
\mathfrak{n}$ are two $g$-orthogonal submodules of the
$A$-eigenspace $\mathfrak{m}_{\lambda}$. Applying Corollary
\ref{CorollaryBrackets}, we conclude that
\[ [\mathfrak{p}_1, [\mathfrak{p}_1,
\mathfrak{p}_2] ] \subset \mathfrak{m}_{\lambda}.
\]
Iterating this
process, we prove that
\[
\mathfrak{n}= [\mathfrak{p}_1, \mathfrak{p}_2]+ [\mathfrak{p}_1,
[\mathfrak{p}_1, \mathfrak{p}_2]]_{\mathfrak{n}}+[\mathfrak{p}_2,
[\mathfrak{p}_1, \mathfrak{p}_2]]_{\mathfrak{n}}+\cdots \subset
\mathfrak{m}_{\lambda}
\]
and $A = \lambda \cdot {\bf 1}$ on $\mathfrak{p}_1 + \mathfrak{p}_2+\mathfrak{n}$.
\end{proof}

\subsection{A  suf\/f\/icient condition for existence of GO-metric}\label{section4.2}

\begin{lemma} Let $M = G/H$ be a  homogeneous space  of a compact Lie group with a reductive decomposition
$\mathfrak{g}= \mathfrak{h} + \mathfrak{m}$.  Assume that  $\Ad_H$-module $\mathfrak{m}$ has a decomposition
\[
\mathfrak{m}= \mathfrak{m}_1+ \cdots + \mathfrak{m}_k
\]
into invariant submodules, such that  for any $i < j$
\[
[\mathfrak{m}_i, \mathfrak{m}_j]=0
\]
or this condition  valid  with one exception $(i,j) =(1,2)$ and  in this case
\[
[\mathfrak{m}_1, \mathfrak{m}_2] \subset \mathfrak{m}_2
\]
and for any $X \in \mathfrak{m}_1$, $Y  \in \mathfrak{m}_2 $ there is $H \in \mathfrak{h}$ such that
$\ad_H Y =  \ad_{X}Y$ and
\[
\ad_H (\mathfrak{m}_1 + \mathfrak{m}_3 + \cdots  + \mathfrak{m}_k)=0.
\]
Then any metric endomorphism  of the form $ A  = \sum x_i \cdot \bold{1}_{\mathfrak{m}_i}$ defines a GO-metric on $M$.
\end{lemma}

\begin{proof} Under the assumptions of lemma,   for $H \in \mathfrak{h} $ and $X_i \in \mathfrak{m}_i$ we have
\begin{gather*}
\Big[H + \sum X_i, \sum x_i X_i\Big]  =  \sum_{i <j}(x_j - x_i)[X_i,X_j] + \sum_i x_i \ad_H X_i
\\
\hphantom{\Big[H + \sum X_i, \sum x_i X_i\Big]}{}
=  (x_2 -x_1 )\ad_{X_1} X_2 + x_2\ad_H X_2+ \ad_H\Bigg(x_1 X_1 + \sum_{k\geq 3} x_k X_k \Bigg).
\end{gather*}
The right-hand side is zero if $H$ is chosen  as in the lemma (where $Y=x_2X_2$ and $X=(x_1-x_2)X_1$).
Now, it suf\/f\/ices to apply Lemma~\ref{GO-criterion}.
\end{proof}

\begin{example} {The homogeneous space  $M = SU_{p+q}/{SU_p \times SU_q}$  is  a  GO-space with respect to any
invariant metric.}
\end{example}

We have the reductive decomposition
\[
\mathfrak{su}_{p+q} =  \mathfrak{h} + \mathfrak{m}= (\mathfrak{su}_p + \mathfrak{su}_q )+ (\mathbb{R}a + \mathfrak{p}),
\]
where $\mathfrak{p} \simeq \mathbb{C}^p \otimes \mathbb{C}^q$ and
$\ad_a |_{\mathfrak{p}} = i \cdot \bold{1}_{\mathfrak{p}}$.
Any metric endomorphism $A = \lambda \cdot \bold{1}_{\mathbb{R}a} + \mu \cdot \bold{1}_{\mathfrak{p}}$
def\/ines a GO-metric by above lemma   since for any $X \in \mathfrak{p}$ there is $H \in \mathfrak{h}$ such that
$\ad_H X = iX$. Note that for $p\neq q$ these manifolds are weakly symmetric spaces~\cite{W1}.

\subsection[GO-metrics on  a  compact group $G$]{GO-metrics on  a  compact group $\boldsymbol{G}$}\label{section4.3}

\begin{proposition} A  compact Lie group $G$ with a left-invariant metric
$g$ is a GO-space if and only if the corresponding Euclidean metric $(\cdot,\cdot)$ on
the Lie algebra $\mathfrak{g}$ is bi-invariant.
\end{proposition}

\begin{proof}
The condition that $(G,g)$ is a GO-space  can be  written as
\[
0=   (X, [X,Y]) = - (\ad_Y X,X) =0.
\]
This shows that the metric $(\cdot,\cdot)$  is bi-invariant.
\end{proof}

Note that a compact Lie group $G$ can admit a non-bi-invariant left-invariant metrics $g$  with homogeneous geodesics.
But  the corresponding GO-space   will have  the form $L/H$  where  the group $L$ will contain $G$ as a proper subgroup.
 See~\cite{DZ} for details.

\section{Homogeneous  GO-spaces  with positive Euler characteristic}\label{section5}

\subsection{Basic facts  about homogeneous manifolds of positive Euler characteristic}\label{section5.1}

Here we recall some properties of homogeneous  spaces with
positive Euler characteristic (see, for example,~\cite{On} or~\cite{BerNik3}).
 A  homogeneous space  $M= G/H$ of a compact connected
  Lie group~$G$ has positive Euler characteristic $\chi(M) >0$ if
  and only if  the  stabilizer $H$ has maximal rank ($\mathrm{rk}(H)
  =\mathrm{rk}(G)$).

 If the  group $G$ acts on $M$  almost ef\/fectively,
 then it is semisimple and the universal covering
   $\widetilde M = \widetilde G/ \widetilde H$
 is a direct product
\[
\widetilde M = G_1/H_1 \times \cdots \times G_k/H_k,
\]
where  $\widetilde G = G_1 \times G_2 \times \cdots \times G_k$ is the decomposition of
the group $\widetilde G$ (which is a covering of~$G$) into a direct product of simple factors
and $H_i = \widetilde H \cap G_i$.

 Any invariant  metric $g^M$ on $M$
  def\/ines an invariant metric $g^{\widetilde M}$ on $\widetilde M$  and the
  homogeneous Riemannian space
  $(\widetilde M = \widetilde G/ \widetilde H, g^{\widetilde M})$ is a direct
  product of homogeneous Riemannian spaces  $(M_i = G_i/H_i, g^{M_i})$, $i=1, \dots, k$,
of simple compact Lie groups $G_i$, see \cite{Kost}.
We have

\begin{proposition}[\cite{Kost}]\label{KostN}
 A  compact almost effective
homogeneous Riemannian space   $(M= G/H,g^M)$ of positive Euler
characteristic is irreducible  if and only if the group~$G$ is
simple. If the group $G$ acts effectively on~$M$, it has trivial center.
\end{proposition}

 This proposition  shows  that  a simply connected  compact GO-space $(M = G/H,
g^M)$ of po\-si\-ti\-ve Euler characteristic is  a direct product of
 simply connected GO-spaces $(M_i = G_i/H_i, g^{M_i})$ of  simple Lie
 groups with positive Euler characteristic. So it is suf\/f\/icient to
 classify  simply connected GO-spaces of a simple  compact Lie group
with positive Euler characteristic.

 A  description of  homogeneous spaces
 $G/H$ of positive Euler
characteristic reduces to description of connected  subgroups $H$ of
maximal rank of  $G$ or equivalently,  subalgebras of maximal rank
of a simple compact Lie  algebra $\mathfrak{g}$, see \cite{BorSieb}
and  also Section~8.10 in \cite{W}.
An important subclass of compact homogeneous spaces of positive
Euler characteristic  consists of  {\bf   f\/lag mani\-folds}. They
are described as adjoint orbits $M = \Ad_Gx$ of a compact connected
semisimple Lie group~$G$ or, in other terms  as quotients  $M=G/H$ of
$G$   by  the centrelizer   $H=Z_G(T)$
 of a~non-trivial torus $T\subset G$.

Note that every compact naturally reductive homogeneous
Riemannian space of positive Euler characteristic is necessarily
normal homogeneous with respect to some transitive semisimple
isometry group \cite{BerNik3}.

\subsection{The main theorem}\label{section5.2}

Let $G$ be a simple compact connected Lie group, $H\subset K\subset
G$ its closed  connected subgroups. We   denote by
 $b= \langle \cdot, \cdot \rangle$ the minus Killing  form
on the Lie algebra $\mathfrak{g}$  and  consider the following $b$-orthogonal decomposition
\[
\mathfrak{g}=\mathfrak{h}\oplus\mathfrak{m}=\mathfrak{h}\oplus\mathfrak{m_1}\oplus\mathfrak{m_2},
\]
where
\[
\mathfrak{k}=\mathfrak{h}\oplus\mathfrak{m_2}
\]
is the Lie algebra of the group $K$. Obviously,
$[\mathfrak{m_2},\mathfrak{m_1}]\subset \mathfrak{m_1}$. Let
$g^M=g_{x_1,x_2}$ be a $G$-invariant Riemannian metric on $M=G/H$,
 generated by the Euclidean metric  $g = (\cdot, \cdot)$ on $\mathfrak{m}$ of the form
\begin{gather}\label{spme1}
g =x_1 \cdot b_{\mathfrak{m}_1}+   x_2 \cdot b_{\mathfrak{m}_2},
\end{gather}
where $x_1$ and $x_2$ are  positive  numbers, or, equivalently, by
the  metric endomorphism
\begin{equation}\label{spme2}
A = x_1 \cdot \bold{1}_{\mathfrak{m}_1} +   x_2 \cdot \bold{1}_{\mathfrak{m}_2}.
\end{equation}

We consider two examples  of such  homogeneous Riemannian spaces $(M= G/H, g_{x_1,x_2})$:
\begin{enumerate}\itemsep=0pt
\item[$a)$] $(G,K,H) = (SO(2n +1),  U(n), SO(2n))$, $n\geq 2$.
The group $G = SO(2 n+1)$  acts transitively on the  symmetric
space $Com(\mathbb{R}^{2 n +2}) = SO(2 n +2)/U(n)$ of
complex structures in $\mathbb{R}^{2 n +2}$ with  stabilizer
$H= U(n) $, see~\cite{Hel}. So we can identify $M = G/H$ with
this symmetric space, but the metric $g_{x_1, x_2}$
is not $SO(2 n +2)$-invariant if $x_2 \neq 2x_1$~\cite{Ker}.

\item[$b)$] $(G,K,H)  = (Sp(n), Sp(1)\cdot Sp(n-1), U(1)\cdot
Sp(n-1))$,  $n\geq 2$.
The group $ G=Sp(n)$ acts transitively on the projective space
$\mathbb{C}P^{2n -1} = SU(2 n + 2)/ U(2 n +1)$ with
stabilizer $H = U(1) \cdot Sp(n -1)$. So we  can identify $M =
G/H$ with $\mathbb{C}P^{2 n -1}$, but the metric $g_{x_1, x_2}$
is not $ SU(2 n + 2)$-invariant if $x_2 \neq 2x_1$, see~\cite{Hel,Ker}.
\end{enumerate}

Now we can state the main theorem  about compact GO-spaces of
positive Euler characteristic.

\begin{theorem}\label{GOPCmain}
Let $(M= G/H,g^M)$ is a simply connected  proper GO-space with positive Euler
characteristic and simple compact  Lie group $G$. Then $M= G/H
=SO(2n+1)/U(n)$, $n\geq 2$, or $G/H=Sp(n)/U(1)\times Sp(n-1)$,
\mbox{$n\geq 2$}, and $g^M = g_{x_1, x_2}$ is any $G$-invariant metric
which is not $G$-normal homogeneous.
The metric $g^M$ is $G$-normal
homogeneous $($respectively, symmetric$)$  when $x_2=x_1$  $($respectively, $x_2=2x_1)$.
Moreover,  these homogeneous spaces  are weakly symmetric flag
manifolds.
\end{theorem}

The non-symmetric  metrics $g_{x_1,x_2}$ have  $G$ as  the full
connected isometry group of the considered  GO-spaces $(M = G/H,
g_{x_1, x_2})$, see discussion
 in \cite{Ker, On, AA}. The claim that all these homogeneous  Riemannian
spaces are weakly symmetric spaces was proved  in \cite{Zi96}.
Note also that Theorem \ref{GOPCmain} allows to simplify some
arguments  in the paper~\cite{BerNik}.

\subsection{Proof of the main theorem}\label{section5.3}

 Using results  from \cite{AA}  and \cite{BerNik}, we reduce  the
 proof to a description of some special decompositions of the root
 system of the Lie algebra $\mathfrak{g}$ of the isometry group~$G$.

Let $M = G/H$ be a homogeneous space of a compact simple Lie group of positive characte\-ris\-tic  and
\[
\mathfrak{g} = \mathfrak{h}+ \mathfrak{m}
\]
associated  reductive decomposition.
The subgroup  $H$ contains a maximal torus  $T$ of $G$.
We consider  the root space decomposition
\[
\mathfrak{g}^{\mathbb{C}} = \mathfrak{t}^{\mathbb{C}} + \sum_{\alpha  \in R} \mathfrak{g}_{\alpha}
\]
of the complexif\/ication $\mathfrak{g}^{\mathbb{C}}$ of the Lie  algebra $ \mathfrak{g}$,
where $\mathfrak{t}^{\mathbb{C}}$
is the Cartan subalgebra associated with $T$  and $R$ is the root system.

For any subset $P \subset R$ we denote by
\[
\mathfrak{g}(P) =
\sum_{\alpha \in P} \mathfrak{g}_{\alpha}
\]
the subspace  spanned
by corresponding root space $\mathfrak{g}_{\alpha}$. Then $H$-module
$\mathfrak{m}^{\mathbb{C}}$ is decomposed into a direct sum
\[
\mathfrak{m}^{\mathbb{C}} = \mathfrak{g}(R_1) + \cdots + \mathfrak{g}( R_k)
\]
of disjoint  submodules, where $R = R_1 \cup \cdots
\cup R_k$ is a disjoint decomposition of $R$  and  subsets~$R_i$ are
symmetric, i.e.\ $-R_i = R_i$. Moreover, real $H$-modules
$\mathfrak{g} \cap \mathfrak{g}(R_i)$ are irreducible. Any invariant
metric on $M$ is def\/ined by  the metric endomorphism $A$ on
$\mathfrak{m}$ whose extension to $ \mathfrak{m}^{\mathbb{C}}$  has
the form
\[
A = \diag( x_1 \cdot \bold{1}_{\mathfrak{p}_1}, \dots , x_{\ell} \cdot \bold{1}_{\mathfrak{p}_{\ell}}),
\]
where $x_i$ are arbitrary positive numbers, $x_i \neq x_j$  and
$\mathfrak{p}_i$ is a direct sum of modules $\mathfrak{g}(R_m)$.

We will assume that $A$ is not a scalar operator (i.e.\ $\ell >1$) and it def\/ines an  invariant metric
with  homogeneous geodesics.  We say that a root  $\alpha$ corresponds to  eigenvalue $x_i$ of $A$
if $ \mathfrak{g}_{\alpha} \subset \mathfrak{p}_i$.

\begin{lemma}  There are two roots $\alpha$, $\beta$ which correspond to different eigenvalues of $A$  such that
$\alpha + \beta$ is a root.
\end{lemma}
\begin{proof} If it is not the case,
$[\mathfrak{p}_1, \mathfrak{p}_i]=0$ for $i \neq 1$
and  $\mathfrak{g}_1 = \mathfrak{p}_1 + [\mathfrak{p}_1, \mathfrak{p}_1]$  would be
a proper ideal of a simple Lie algebra $\mathfrak{g}$.
\end{proof}

Now, consider the roots $\alpha$ and $\beta$ as in the previous lemma.
Since $R(\alpha, \beta) := R \cap \spann \{ \alpha, \beta \}$
is  a rank 2  root  system,  we can always  choose  roots $\alpha, \beta \in R$ which  form a basis of the root
system $R(\alpha, \beta)$.  Then the subalgebra
\[
\mathfrak{g}_{\alpha, \beta} := \mathfrak{t}^{\mathbb{C}} + \sum\limits_{\gamma \in
R(\alpha,\beta)} \mathfrak{\mathfrak{g}_{\gamma}}
\]
of $\mathfrak{g}^{\mathbb{C}}$ is the centralizer of  the  subalgebra
$\mathfrak{t}' = \ker \alpha \cap \ker \beta \subset \mathfrak{t}^{\mathbb{C}}$.

Then the orbit $G_{\alpha, \beta}o  \subset M$ of the corresponding
subgroup $G_{\alpha, \beta}= T' \cdot G'_{\alpha, \beta} \subset G$
is a totally geodesic  submanifold (see Corollary \ref{togc1}),
hence a proper GO-space  with the ef\/fective  action of the  rank two
simple group $G'_{\alpha, \beta} $ associated with the root system
$R(\alpha, \beta)$ (see Proposition \ref{tot}). Note that it has
positive Euler characteristic since the  stabilizer of  the point
$o$ contains the two-dimensional torus generated by vectors
$H_{\alpha}, H_{\beta} \in \mathfrak{t}^{\mathbb{C}}$ associated
with roots $\alpha$, $\beta$. Recall
that
$H_{\alpha } = \frac{2}{\langle \alpha, \alpha\rangle}\,  b^{-1} \cdot \alpha$.

\begin{proposition}\label{spr2}
Every proper GO-space $(M = G/H, g^M)$  with positive Euler characteristic of a simple group $G$ of
rank~$2$  is locally isometric to the manifold $M =
SO(5)/U(2)$  with the metric defined  by  the metric endomorphism
\[
A = x_1 \cdot \bold{1}_{\mathfrak{g}(R^{s})}  +   x_2 \cdot
\bold{1}_{\mathfrak{g}(R^{\ell})} , \qquad x_1 \neq x_2 >0
\]
where
\[
R^{s}= \{ \pm \epsilon_1, \pm \epsilon_2\},\qquad  R^{\ell} = \{ \pm
\epsilon_1 \pm \epsilon_2 \},
\]
are the sets of short and, respectively, long roots of the Lie algebra $\mathfrak{so}(5)$.
We  may assume  also that
\[
\mathfrak{m}^{\mathbb{C}} = \mathfrak{g}(R^{s} \cup \{ \epsilon_1 + \epsilon_2\})\qquad \mbox{and}\qquad
\mathfrak{h}^{\mathbb{C} } = \mathfrak{t}^{\mathbb{C}} + \mathfrak{g}_{\epsilon_1 - \epsilon_2}.
\]
\end{proposition}

\begin{proof}
Proof of this proposition follows from results of the papers
\cite{AA} and \cite{BerNik}. Indeed, the group $G$  has the Lie algebra  $\mathfrak{g}$
isomorphic to  $su(3)=A_2$,
$so(5)=sp(2)=B_2=C_2$ or  $g_2$. Since the universal Riemannian
covering of a GO-space  is a GO-space (Lemma \ref{GO-criterion}), we may assume
without loss of generality that $G/H$ is simply connected.

If $\mathfrak{g}=su(3)$, then $G/H=SU(3)/S(U(2)\times U(1))$ (a
symmetric space) or $G/H =SU(3)/T^2$, where $T^2$ is a maximal torus
in $SU(3)$. Both these spaces are f\/lag manifolds, and  results
of \cite{AA}  show that any GO-metric on these spaces is
$SU(3)$-normal homogeneous.

If $\mathfrak{g}=so(5)=sp(2)$, then
$(\mathfrak{g},\mathfrak{h})=(so(5), \mathbb{R}^2)$,
$(\mathfrak{g},\mathfrak{h})=(so(5),\mathbb{R}\oplus su(2)_l)$,
$(\mathfrak{g},\mathfrak{h})=(so(5),\mathbb{R}\oplus su(2)_s)$,
or $(\mathfrak{g},\mathfrak{h})=(so(5),su(2)_l \oplus su(2)_l)$,
where $su(2)_l$,  (respectively, $su(2)_s$) stands for a~three-dimensional subalgebras generated by all
long (respectively, short) roots of $\mathfrak{g}$.
The last pair corresponds to the irreducible symmetric space $SO(5)/SO(4)$, which admits
no non-normal invariant metric. All other spaces are f\/lag manifolds. Results of~\cite{AA} implies that the only possible pair is
$(\mathfrak{g},\mathfrak{h})=(so(5),\mathbb{R}\oplus su(2)_l)$, which corresponds
to the space $SO(5)/U(2)=Sp(2)/U(1)\cdot Sp(1)$.

For $\mathfrak{g}=g_2$ the statement of proposition is proved  in
\cite[Proposition 23]{BerNik}.
\end{proof}

\begin{corollary}\label{maint}
Let $G$ be a simple compact Lie group   and  $M=G/H$ a  proper GO-space  with positive Euler characteristic.
 Then the root system $R$ of the complex Lie algebra
$\mathfrak{g}^{\mathbb{C}}$  admits a~disjoint decomposition
\[
R = R_0 \cup R_1 \cup R_2,
\]
where $R_0$ is the root system of the  complexified stability
subalgebra $ \mathfrak{h}^{\mathbb{C}} $,with the following
properties:
\begin{enumerate}
\itemsep=0pt
\item[$i)$]  If $\alpha \in R_1$, $\beta \in R_2$  and $ \alpha + \beta  \in
R$ then $\alpha - \beta  \in R$  and the  rank~$2$ root system $R(\alpha,
\beta)$ has type $B_2 = C_2$.

\item[$ii)$] Moreover, if $\alpha$, $\beta$ is a basis of $R(\alpha, \beta)$
(that is $\langle \alpha, \beta \rangle <0$), then one of the roots
$\alpha$,  $\beta$ is short  and the  other is long  and  one of the  long
roots $ \alpha \pm \beta$ belongs to
$R_0$  and second one belongs to $ R_1 \cup R_2$.

\item[$iii)$] If both roots $\alpha$, $\beta$ are short, then  one  of the
long roots $\alpha \pm \beta$ belongs to $ R_0$ and the  other  belongs
to
$R_1 \cup R_2$.

\item[$iv)$]  If $\alpha \in R_1$ and $\beta \in R_2$  are long roots, then
$\alpha \pm \beta \notin R$.
\end{enumerate}
\end{corollary}

We will call  a  decomposition with the above properties   a  {\bf special decomposition}.
 Corollary~\ref{maint} implies

\begin{corollary}\label{olr}
There is no proper GO-spaces of positive Euler characteristic  with
simple isometry group $G = SU(n),\, SO(2n), \,E_6,\,E_7,\,E_8$ $($these are
all simple Lie algebras with all roots of the  same length   $($simply-laced root system$))$.
\end{corollary}

\begin{corollary}[\protect{\cite[Proposition 23]{BerNik}}]\label{g2}
Any GO-space $(G/H, \mu)$ of positive Euler characteristic with $G=G_2$
is normal homogeneous.
\end{corollary}

Now, we describe   all { \bf special decompositions} of the  root
systems of types $B_n$, $C_n$, $F_4$. We will use notation from~\cite{GOV} for root  systems  and simple roots.

\begin{lemma}
 The root system
\[
R(F_4) = \{ \pm \epsilon_i, \ 1/2( \pm \epsilon_1 + \mp \epsilon_2 \pm \epsilon_3 +
\pm \epsilon_4, \pm \epsilon_i \pm \epsilon_j ),\  i,j =1,2,3,4,\  i \neq j \}
\]
does not admit a special decomposition.
\end{lemma}

\begin{proof} Assume that  such a decomposition exists. Then we can choose  roots $\alpha \in R_1$, $\beta \in R_2$
such that
$\alpha  \pm \beta $ is a root.  Then $\alpha$, $\beta$ has  dif\/ferent length  and we may assume that
 $|\alpha| < |\beta|$ and  $ \langle \alpha, \beta \rangle <0$.
 Then we can include $\alpha$, $\beta$ into a system of simple roots
$ \delta$, $\alpha$, $\beta$, $\gamma $, see~\cite{GOV}.  Since  all  such systems  are conjugated, we may assume that
 $\alpha = \epsilon_4$, $\beta = -\epsilon_4 + \epsilon_3 $, see \cite{GOV}. Then  we get contradiction, since
 $\alpha - \beta $ is not a root.
\end{proof}

Now we describe two special decompositions  for  the root systems
\[
R(B_n) = \{ \pm \epsilon_i, \ \pm \epsilon_i \pm \epsilon_j, \  i,j = 1, \dots, n   \}
\]
and
\[
R(C_n) = \{ \pm 2\epsilon_i,  \ \pm \epsilon_i \pm \epsilon_j,  \  i,j = 1, \dots ,n  \}
\]
of types $B_n$ and $C_n$. Note  that in  both  cases $R_A = \{ \pm(\epsilon_i - \epsilon_j)  \}$
is a closed subsystem. We set  $R_A^{+} = \{ \pm(\epsilon_i + \epsilon_j) \}$.

We  denote by $R^+$ the  standard  subsystem of positive roots  of a
root system $R$  and by~$R^s$ and~$R^{\ell}$  the subset of  short
and, respectively, long
 roots of $R$. Then   there is  a special  decomposition
   $R = R_0 \cup R_1 \cup R_2$
of the  systems $R(B_n)$, $R(C_n)$  which we call  the  standard
decomposition:
\begin{gather*}
R(B_n)  = R_A \cup R^{s}\cup R_A^{+},
\\
R(C_n)  = R_A \cup R^{\ell}\cup R_A^{+}.
\end{gather*}
These  decompositions  def\/ine  the  following  reductive decompositions of the homogeneous spaces
$SO(2n+1)/U(n)$  and $Sp(n)/ U(n)$:
\begin{gather*}
\mathfrak{so}(2n+1) = \mathfrak{h} + (\mathfrak{m_1}+ \mathfrak{m}_2)=
\mathfrak{g}(R_A) + (\mathfrak{g}(R^s) + \mathfrak{g}(R_A^+)),
\\
 \mathfrak{sp}(n) = \mathfrak{h} + (\mathfrak{m_1}+ \mathfrak{m}_2)=
\mathfrak{g}(R_A) + (\mathfrak{g}(R^{\ell}) + \mathfrak{g}(R_A^+)),
\end{gather*}
where $\mathfrak{m}_1$, $\mathfrak{m}_2$  are irreducible
submodules of $\mathfrak{m}$. It is known \cite{AA} that any  metric
endomorphism $A = \diag( x_1 \cdot \bold{1}_{\mathfrak{m}_1},  x_2
\cdot \bold{1}_{\mathfrak{m}_2} )$ def\/ines a metric with homogeneous
geodesics on the  corresponding manifold $M = G/H$ (see a discussion
before the statement of Theorem~\ref{GOPCmain}). Now, the proof of
Theorem~\ref{GOPCmain} follows from the following proposition.

\begin{proposition}\label{spdec}
Any special decomposition of the root systems $R_B$, $R_C$ is conjugated to the standard one.
\end{proposition}

\begin{proof}  We give a proof of this proposition for $R(B_n)$. The  proof for $R(C_n)$ is similar.

 Let
\[  R(B_n)  = R_0 \cup R_1 \cup R_2
\]
be a special decomposition of $R(B_n)$. We may assume that there are roots $\alpha \in R_1$ and $\beta \in R_2$ with
$\langle \alpha, \beta \rangle <0$ and $|\alpha| < |\beta|$.
Then we can include $\alpha, \beta$ into a system of simple  roots,
which, without loss of
generality,  can be written as
\[
\alpha_1 = \epsilon_1 - \epsilon_2, \ \dots , \ \epsilon_{n-2} - \epsilon_{n-1}, \
 \epsilon_{n-1} + \epsilon_n = \beta , \ -\epsilon_n = \alpha.
\]
Then $ (\epsilon_{n-1} - \epsilon_n) \in R_0$.
We need the following lemma.

\begin{lemma} Let $  R(B_n)  = R_0 \cup R_1 \cup R_2$   be a special decomposition as  above,
$V'= \epsilon_n^\perp$  the orthogonal complement  of  the vector $\epsilon_n$  and
$ R(B_{n-1})= R': = R \cap V' $    the   root  system  induced in  the hyperspace  $V'$. Then   the
induced decomposition
$  R' = R'_0 \cup R'_1 \cup R'_2 $,  where  $ R'_i := R_i \cap V'$,  is a special  decomposition.
\end{lemma}

\begin{proof} It is suf\/f\/icient to check that  subsets  $R'_1$, $R'_2$ are not empty.

We say that  two roots $\gamma$, $\delta$ are  $R_0$-{\bf equivalent} ($\gamma \sim \delta$) if their
dif\/ference belongs to $R_0$. The equivalent roots belong to the same component
$R_i$. The  root  $\epsilon_{n-1} = \epsilon_n - (\epsilon_{n-1} - \epsilon_n)$  is $R_0$-equivalent to
$\alpha = \epsilon_n$. Hence it belongs to $R_1$.

We say that a pair of  roots $\gamma$, $\delta$  with $\langle \gamma, \delta \rangle <0$ is
{\bf  special} if   one of  the roots  belongs to $R_1$ and another to~$R_2$. Then  they have dif\/ferent length (say, $|\gamma| < |\delta|$). Moreover, the  root $\gamma + \delta$ is
short  and it belongs to the same  part $R_i$, $i =1,2$ as  the short root $\delta$  and  the root  $2\gamma + \delta $
is  long  and it belongs to  $R_0$.

Consider the  roots $\sigma_{\pm}= \pm \epsilon_{n-2}+  \epsilon_{n-1}$. They have negative scalar product with
$\epsilon_{n-1} \in R_1$ and $\beta = \epsilon_{n-1} +  \epsilon_{n} \in R_2$. They can not belong to $R_1$
since then we get  a special pair $\delta_{\pm}, \beta $ which consists
of  long roots.  They  both can not belong to $R_0$ since otherwise  the  root
$\epsilon_{n-2}  \sim \epsilon_{n-1}  \in R_1$ and
$ \pm \epsilon_{n-2} + \epsilon_n \sim  \epsilon_{n-1} + \epsilon_n \in R_2$ and we get a special pair
\[
\gamma = \epsilon_{n-2} \in R_1, \qquad  \delta = -\epsilon_{n-2}+ \epsilon_n \in R_2,
\]
such that $2 \gamma + \delta \in R_0$, which is impossible.
We conclude that one of the  roots $ \sigma_{\pm}= \pm \epsilon_{n-2} + \epsilon_{n-1} \in R' $ must belongs to $R_2$.
Since the root  $\epsilon_{n-1} \in R'$ belongs to $R_1$, the lemma is proved.
\end{proof}

Now we prove  the proposition by induction  on $n$.  The claim is true for  $n=2$ by Proposition~\ref{spr2}.
Assume that it is true for $R(B(n-1))$  and let  $  R(B_n)  = R_0 \cup R_1 \cup R_2$ be a special decomposition as above.
 By lemma, the  decomposition $  R' = R'_0 \cup R'_1 \cup R'_2 $,  indiced in the hyperplane $V' = e_n\perp,$ is  a special
decomposition.
By  inductive hypothesis we may assume that it  has the standard form:
\[
R_0 = \{\pm( \epsilon_i - \epsilon_j) \},  \qquad  R_1 = \{ \pm \epsilon_i \},\qquad
R_2 = \{ \pm(\epsilon_i + \epsilon_j) ,\  i,j =1, \dots, n-1\}.
\]
This implies  that  the initial decomposition is also standard.
\end{proof}

\subsection*{Acknowledgements}
The f\/irst author was partially supported by the Royal Society
(Travel Grant 2007/R3). The second author was partially supported
by the State Maintenance Program for the Leading Scientif\/ic
Schools of
the Russian Federation (grant NSH-5682.2008.1).
We are grateful to all referees, whose comments and suggestions
permit us to improve the presentation of this article.

\pdfbookmark[1]{References}{ref}
\LastPageEnding

\end{document}